\documentclass[12pt, twoside]{article}
\usepackage{amsmath,amsthm,amssymb}
\usepackage{times}
\usepackage{enumerate}

\def\subjclass#1{{\renewcommand{\thefootnote}{}%
\footnote{\hspace{-0.6cm}\emph{Mathematics Subject Classification
(2020):} #1}}} 
\def\subj#1{{\renewcommand{\thefootnote}{}%
\footnote{\hspace{-0.6cm}\emph{keywords:} #1}}}

\newtheorem{thm}{Theorem}[section]
\newtheorem{proposition}[thm]{Proposition}
\newtheorem{theorem}[thm]{Theorem}
\newtheorem{corollary}[thm]{Corollary}
\newtheorem{lemma}[thm]{Lemma}

\theoremstyle{definition}
\newtheorem{definition}[thm]{Definition}
\newtheorem{remark}[thm]{Remark}
\newtheorem{example}[thm]{Example}

\numberwithin{equation}{section}

\begin{document}


\baselineskip=17pt


\title{\bf\textsc{Some properties of the generalized $p$-adic gamma function}}

\author{Rafik BELHADEF and Nour Elhouda SAHALI}

\date{}

\maketitle

 \subj{$p$-adic number, $p$-adic factorial, $p$-adic gamma function}


\begin{abstract}

	In this paper, we define a q-adic factorial and we demonstrate some properties of a generalized $p$-adic gamma function. Also, some numerical examples have been given.

\end{abstract}


\section{Introduction}

The $p$-adic gamma function $\Gamma_{p}$ is a $p$-adic integer function analogous to the classical gamma function. In 1975, Morita defined $\Gamma_{p}$ explicitly by:

\begin{equation*}
	\Gamma_{p}(n)=(-1)^n\prod_{j=1,(p,j)=1}^{n-1}j
\end{equation*}

There were several generalizations of the $p$-adic gamma function (see \cite{ma}, \cite{8} and \cite{ko}), one of which was introduced by Kaori Ota \cite{ko} in 1994. 
To study the generalized hypergeometric function, Ota defined the generalized $p$-adic gamma function $\Gamma_{q}$ with $q=p^t, t\in \mathbb{N^*}$ by the formula:

\begin{align*}
	\Gamma_{q}(x+1)=\prod_{\ell=0}^{t-1} \Gamma_{p}(h_{\ell}(x)+1) 
\end{align*}

such that $ x=\sum_{j=0}^{\infty}  x_{j} p^{j}$ and $h_{\ell}(x)= \sum_{j \geq \ell} x_{j} p^{j-\ell}$, for $\ell \in \mathbb{N}$.

Another extension for the $p$-adic gamma function was adopted by N. Koblitz , He use the same notation $\Gamma_{q}$ for a different function, but this is not confusing. Both Koblitz and Ota have defined without giving all the properties of the generalized $p$-adic gamma function.\\

In our paper, we first demonstrate some propositions given by Ota
(see  Proposition \ref{pr2}, Proposition \ref{pr.4}, and Proposition \ref{pr3}). 
Secondly, we define q-adic factorial and we use this concept to demonstrate the combinatorial properties of the generalized $p$-adic gamma function, similar to that of $p$-adic gamma function (see  Proposition \ref{pr4.2.2}, Proposition \ref{pr1}, and Corollary \ref{corl4.2.5}). 
Next, we propose an expose a Mahler expansion of $\Gamma_{q}$, and we prove the relationship between its coefficients (see Proposition \ref{4.2.8}). Finally, some numerical examples are given (Examples \ref{exa22}).

\section{Preliminary}

Throughout this paper We use the following concepts: $p$ is a prime number, $\mathbb{Z}$ is the set of all the real integers, $\mathbb{Z_{-}}$ (resp. $\mathbb{Z_{+}}$) is the set of all the negative real integers (resp. all the positive real integers), $\mathbb{N}$ is the set of all the non-negative integers, $\mathbb{Q}$ is the field of rational numbers, and $\mathbb{R}$ is the field of real numbers. We use $\left|.\right|$ to denote the absolute value in $\mathbb{R}$, $[.]$ the real integer part, $\nu_p$ the $p$-adic valuation, and $\left|.\right|_p$ the $p$-adic absolute value. The field of $p$-adic numbers $\mathbb{Q}_{p}$ is the completion of $\mathbb{Q}$ with respect to the $p$-adic absolute
value. The ring of $p$-adic integers $\mathbb{Z}_{p}$  contains the p-adic numbers which satisfy $\left|x\right|_p \leq 1$.


\subsection{$p$-adic Factorial and $p$-adic Gamma Function}

In this subsection, we define the $p$-adic factorial function, $p$-adic gamma function, and some of their basic properties, to be needed in the next section.

\begin{definition} \label{Analog-fact}\cite{rb}
	The $p$-adic factorial is defined by $ 0!_p=1$ and for $n \in \mathbb{N^*}$ by
	\begin{equation}
		n!_{p}=\underset{\underset{(p,j)=1}{j=1}}{\overset{n}{\prod }} j 
	\end{equation}
\end{definition}

The p-adic gamma function has been well used in dynamic systems and string theory. This function is studied by \cite{8} and \cite{10}, to give some properties of polynomials.\\
The function $ n!$ cannot be extended by continuity on $\mathbb{Z}_p$, because $\underset{n\rightarrow +\infty }{\lim }n!=0$ in $\mathbb{Z}_p$. So, we have the definition of $p$-adic gamma function as follows:

\begin{definition} \label{Analog-gamma}\cite{11}
	The $p$-adic gamma function is defined by Morita as the continuous function 
	\begin{equation*}
		\Gamma _{p}:
		\mathbb{Z}
		_{p}\longrightarrow 
		\mathbb{Z}_{p}
	\end{equation*}
	as an extension of the following sequence, with $n \in \mathbb{Z_{+}}$
	\begin{equation}
		\Gamma_{p}(n)=(-1)^n\prod_{j=1,(p,j)=1}^{n-1}j
	\end{equation} 
	Furthermore,
	$$\Gamma _{p}(z)=\underset{n\underset{\text{in }\mathbb{Z}_{p}}{\rightarrow }z}{\lim }\Gamma _{p}(n)
	=\underset{n\underset{\text{in } \mathbb{Z}_{p}}{\rightarrow }z}{\lim }\left( -1\right) ^{n}\underset{\underset{(p,j)=1}
		{j=1}}{\overset{n-1}{\prod }}j$$
\end{definition}

Here, we cite some properties of $\Gamma_{p}$ that we need to prove the theorems in the next section.

\begin{proposition} \label{prop1-gamma} \cite{1}
	The function $\Gamma_{p}$ satisfies the following properties:
	\begin{enumerate}
		\item[1)] $\Gamma_{p}(0)=1$ , $\Gamma_{p}(1)=-1$, $\Gamma_{p}(2)=1$
		\item[2)] $\Gamma_{p}(n+1)=(-1)^{n+1} n!_{p}$, $\forall n \in \mathbb{N}$
	\end{enumerate}
\end{proposition}

Other important arithmetic formulas are given in the following proposition:

\begin{proposition} \label{prop5-gamma} \cite{1}
	Let $n \geq 1$, its $p$-adic expansion be $\displaystyle\sum_{i=0}^{\ell }n_{i}p^{i}$, and the sum of digits be $S_n=\displaystyle\sum_{i=0}^{\ell}n_{i}$. Then,
	\begin{enumerate}
		\item[1)] $\Gamma _{p}(n+1)=\dfrac{\left( -1\right) ^{n+1}n!}{\left[ \frac{n}{p}\right]! \times p^{\left[ \frac{n}{p}\right] }}$. In particular, $\Gamma _{p}(p^{n})=\dfrac{\left( -1\right) ^{p}p^{n}!}{p^{n-1}! \times p^{p^{n-1}}}$.
		\item[2)] $\Gamma _{p}(np+k+1)=\dfrac{\left( -1\right) ^{np+k+1}\left( np+k\right) !}{%
			n! \times p^{n}}$ , \ \ for  \ \  $0\leq k<p$.
		\item[3)] \bigskip $n!=\left( -1\right) ^{n+1-\ell}\left( -p\right) ^{\frac{n-S_{n}}{p-1}}%
		\underset{i=0}{\overset{\ell}{\prod }}\Gamma _{p}\left( \left[ \frac{n}{p^{i}}%
		\right] +1\right) $.
	\end{enumerate}
\end{proposition}

\begin{theorem} \label{theorem6-gamma} \cite{1} \textbf{(Mahler expansion of $\Gamma_{p}$)}
	
	Let the Mahler expansion of $\Gamma_{p}$:
	$$\Gamma_{p}(x+1)=\sum_{k \geq 0} \alpha_{k} \binom{x}{k}$$
	The coefficients $\alpha_{k}$ verify the following relationship:
	\begin{align}
		\exp\left(x+\frac{x^{p}}{p}\right) \frac{1-x^{p}}{1-x}=\sum_{k \geq 0} (-1)^{k+1} \alpha_{k} \frac{x^{k}}{k!} 
	\end{align}
\end{theorem}


\section{Main Results and Proofs}

\bigskip

Inspired by the works of Belhadef \cite{rb}, we will establish a definition of a $q$-adic factorial. We also define the generalization of the $p$-adic gamma function and we demonstrate their properties.
Throughout this section, we consider $q=p^{t}$, for a positive integer $t$. 

\subsection{$q$-adic factorial}

\begin{definition} \label{q-Analog-fact}
	The $q$-adic factorial is defined by $0!_q=1$ and for $ n> 0$ by
	\begin{equation}
		n!_{q}=\underset{\underset{(q,j)=1}{j=1}}{\overset{n}{\prod }} j 
	\end{equation}
\end{definition}

\begin{remark} 
	If $1 \leq n \leq q-1$, then $(q,j)=1$, for all $1 \leq j \leq n$. So, $n!_{q}=n! $.
\end{remark}

\begin{lemma}\label{two-adic-fact}
	We have $(qk)!_{q}=(qk-1)!_{q} .$
\end{lemma}
\begin{proof} 
	We observe that $(q,qk)\neq 1$, it means that $$(qk)!_{q}=\underset{\underset{(q,j)=1}{j=1}}{\overset{qk}{\prod }} j=\underset{\underset{(q,j)=1}{j=1}}{\overset{qk-1}{\prod }} j=(qk-1)!_{q} .$$
\end{proof}

\begin{example} \label{exa22} 
	In Tables 1, 2, and 3, we calculate some $q$-adic factorials of some positive integers. For $q=2^{2}, 2^{3}, 3^{2}$.	

$$\textbf{Table 1: The $2^{2}$-adic factorial }$$
\begin{equation*}
    \begin{tabular}{ccccccccccccc}
    \hline   \textbf{$n$} & \textbf{0} & \textbf{1} & \textbf{2} & \textbf{3} & \textbf{4} & \textbf{5} & \textbf{6} & \textbf{7} & \textbf{8} & \textbf{9} & \textbf{10} & \textbf{11} \\
    \hline  \textbf{$n!_{2^{2}}$} & 1 & 1 & 2 & 6 & 6 & 30 & 180 &  1260  &  1260   & 11340  &  113400 & 1247400 \\
		\hline
    \end{tabular}
\end{equation*} 

$$\textbf{Table 2: The $2^{3}$-adic factorial }$$
\begin{equation*} 
    \begin{tabular}{ccccccccccccc}
		\hline \textbf{$n$} & \textbf{0} & \textbf{1} & \textbf{2} & \textbf{3} & \textbf{4} & \textbf{5} & \textbf{6} & \textbf{7} & \textbf{8} & \textbf{9} & \textbf{10} & \textbf{11}\\
		\hline	\textbf{$n!_{2^{3}}$} & 1 & 1 & 2 & 6 & 24 & 120 & 720 &  5040  &  5040   & 45360  &  453600 & 4989600 \\
		\hline
    \end{tabular} 
\end{equation*} 

$$\textbf{Table 3: The $3^{2}$-adic factorial }$$
\begin{equation*} 
    \begin{tabular}{ccccccccccccc}
\hline \textbf{$n$} & \textbf{0} & \textbf{1} & \textbf{2} & \textbf{3} & \textbf{4} & \textbf{5} & \textbf{6} & \textbf{7} & \textbf{8} & \textbf{9} & \textbf{10} & \textbf{11}\\
		\hline	\textbf{$n!_{3^{2}}$} & 1 & 1 & 2 & 6 & 24 & 120 & 720 &  5040  &  40320   & 40320  &  403200 & 4435200 \\
		\hline
    \end{tabular} 
\end{equation*} 
\end{example}

\bigskip

The next theorem represents a $q$-generalization of the Wilson congruence, which is the key to some results in this section.

\begin{theorem} \label{wilson-cong2}
	Let $a \in \mathbb{Z}$ and $s \in \mathbb{Z_{+}}$. Then
	\begin{enumerate}
		\item[1)] For $p\geq 3$ and $s \geq 1$, we have $\underset{\underset{(q,j)=1}{j=a}}{\overset{a+p^s-1}{\prod }} j \equiv -1\left( {mod}\text{ }p^{s}\right)$. \label{3.3}
		\item[2)] For $p=2$ and $s \geq 3$, we have $\underset{\underset{(q,j)=1}{j=a}}{\overset{a+2^s-1}{\prod }} j \equiv 1\left( {mod}\text{ }2^{s}\right).$
	\end{enumerate}
\end{theorem}		

From this $q$-generalization, we obtain the following congruence:

\begin{proposition} \label{cong-fact2}
	Let $n \in \mathbb{N}$ and $s \in \mathbb{Z_{+}}$.
	\begin{enumerate}
		\item[1)] If $p\geq 3$ and $s \geq 1$, then $\dfrac{(n+p^s)!_{q}}{n!_{q}} \equiv -1\left( {mod}\text{ }p^{s}\right)$.
		\item[2)] If $p=2$ and $s \geq 3$, then $\dfrac{(n+2^s)!_{q}}{n!_{q}}\equiv 1\left( {mod}\text{ }2^{s}\right)$.
	\end{enumerate}
\end{proposition}
\begin{proof} 
	We have $$\dfrac{(n+p^s)!_{q}}{n!_{q}}=\underset{\underset{(q,j)=1}{j=n+1}}{\overset{n+p^s}{\prod }} j $$
	From case 1 of Theorem \ref{wilson-cong2} with $a=n+1$, we obtain the congruence for $p\geq 3$ and $s \geq 1$. From case 2 of the same Theorem with $a=n+1$, we obtain the congruence for $p=2$ and $s \geq 3$.
\end{proof}

More generally, we have the following theorem and its immediate corollaries:

\begin{proposition} \label{cong-fact-gen2}
	Let $n \in \mathbb{N}$, and $m, s \in \mathbb{Z_{+}}$.
	\begin{enumerate}
		\item[1)] If $p\geq 3$ and $s \geq 1$, then $\dfrac{(n+mp^s)!_{q}}{n!_{q}}  \equiv (-1)^m\left( {mod}\text{ }p^{s}\right)$.
		\item[2)] If $p=2$ and $s \geq 3$, then $\dfrac{(n+m2^s)!_{q}}{n!_{q}}\equiv 1\left( {mod}\text{ }2^{s}\right)$.
	\end{enumerate}
\end{proposition}

\begin{proof}
	The proof is done by induction on $m$.
\end{proof}

\begin{corollary} \label{norm3-fact2}
	For $p\geq 3$, $n \in \mathbb{N}$ and $s \in \mathbb{Z_{+}}$, we have $|n!_{q}|_{p}=1$ and
	$$|(n+p^s)!_{q}+n!_{q}|_{p}\leq \frac{1}{p^s}$$
\end{corollary}

\begin{corollary}  \label{norm2-fact2} 
	For $p=2$, $n \in \mathbb{N}$ and $s \in \mathbb{Z_{+}}$ with $s \geq 3$, we have $|n!_{q}|_{2}=1$ and
	$$|(n+2^s)!_{q}-n!_{q}|_{2}\leq \frac{1}{2^s}$$
\end{corollary}


\subsection{Generalized $p$-adic Gamma Function}

In this part, we introduce the definition of the generalized $p$-adic gamma function and some properties given by Kaori Ota in \cite{ko} and we demonstrate other properties of this function that do not exist in the paper of Ota.

\begin{definition}\label{def1} \cite{ko}
	Let $x \in \mathbb{Z}_{p}$ given by a $p$-adic expansion $x=\sum_{j=0}^{\infty}  x_{j} p^{j}$, $x_{j} \in \{0,1,\cdots,(p-1)\}$, $x_{0} \neq 0$. For $\ell \in \mathbb{N}$, a map $h_{\ell}: \mathbb{Z}_{p} \longrightarrow \mathbb{Z}_{p}$ is define by:
	\begin{align}
		h_{\ell}(x)= \sum_{j \geq \ell} x_{j} p^{j-\ell}  \label{3.25}
	\end{align}
	
	in this case, we have 
	\begin{align}
		x &= \sum_{j=0}^{\ell-1} x_{j}p^{j}+p^{\ell}h_{\ell}(x). \label{2.23}
	\end{align}
\end{definition}

\begin{definition}\cite{ko}\label{def2}
	For $x \in \mathbb{Z}_{p}$, the generalized $p$-adic gamma function is defined from $\mathbb{Z}_{p}$ to $\mathbb{Z}_{p}^{*}$ by:
	\begin{align}
		\Gamma_{q}(x+1)=\prod_{\ell=0}^{t-1} \Gamma_{p}(h_{\ell}(x)+1) \label{3.24}
	\end{align}
\end{definition}


\begin{remark}
	\label{rm1} \noindent
	For $t=1$ the function $\Gamma _{q}$ coincides with $\Gamma _{p}$,  i.e.
	\begin{equation*}
		\Gamma _{q}(x+1)=\Gamma _{p}(h_{0}(x)+1)=\Gamma _{p}(x+1),
	\end{equation*}
	
\end{remark}

\begin{example}
	\label{ex411} For $t=2$, we have $q=p^{2}$ so 
	\begin{equation*}
		\Gamma _{q}(x+1)=\Gamma _{p}(h_{0}(x)+1)\Gamma _{p}(h_{1}(x)+1)=\Gamma
		_{p}(x+1)\Gamma _{p}(h_{1}(x)+1),
	\end{equation*}%
	according to the relationship (\ref{2.23}) we have $x=x_{0}+ph_{1}(x)$ hence $h_{1}(x)=\frac{x-x_{0}}{p}$. Therefore 
	\begin{equation*}
		\Gamma _{q}(x+1)=\Gamma _{p}(x+1)\Gamma _{p}\left( \frac{x-x_{0}}{p}%
		+1\right).
	\end{equation*}%
	For example, for $x=3$, $p=3$, $t=2$, $q=3^{2}$. Then, $x_{0}=1$ so 
	\begin{equation*}
		\Gamma _{9}(4)=\Gamma _{3}(4)\Gamma _{3}(2)=2
	\end{equation*}%
\end{example}

\begin{proposition}\cite{ko} \label{pr2} 
	\begin{enumerate}
		\item For a positive integer $n$, we have 
		\begin{equation}
			\Gamma _{q}(n+1)=(-1)^{A_{n}}p^{-B_{n}}n!_{q},
		\end{equation}%
		where 
		\begin{equation*}
			A_{n}=t+\sum_{i=0}^{t-1}t\left[ \frac{n}{p^{i}}\right] \; and \;B_{n}=%
			\sum_{i=1}^{t}\left[ \frac{n}{p^{i}}\right] -t\left[ \frac{n}{p^{t}}\right] 
		\end{equation*}%
		\label{pr.2}
		\item We have the functional equation
		\begin{equation}
			\Gamma _{q}(x+1)=%
			\begin{cases}
				(-1)^{k+1}h_{k}(x)\Gamma _{q}(x) & if \;x\notin q\mathbb{Z}%
				_{p}\;et\;k=v_{p}(x), \\ 
				(-1)^{t}\Gamma _{q}(x) & if i\;x\in q\mathbb{Z}_{p}.%
			\end{cases}%
		\end{equation}
	\end{enumerate}
\end{proposition}

\begin{proposition}
	\textbf{(The Complement formulas)} \cite{ko} \label{pr.4} 
	\begin{enumerate}
		\item For $p\neq 2$ and $\forall x\in \mathbb{Z}_{p}$ , we have
		\begin{equation}
			\Gamma _{q}(x)\Gamma _{q}(1-x)=(-1)^{t-1+R_{t}(x)}, \label{3.34}
		\end{equation}%
		where, $R_{t}(x)\in \{1,2,\cdots ,q\}$, such as $R_{t}(x)\equiv x(\mod q)$.%
		\item For $p=2$ 
		\begin{equation}
			\Gamma _{q}(x)\Gamma _{q}(1-x)=%
			\begin{cases}
				(-1)^{t+1+x_{t}} & if \;x\notin q\mathbb{Z}_{2}\;et\;x\notin \mathbb{Z}%
				_{2}^{\ast }, \\ 
				(-1)^{t+x_{t}} & if \;x\in q\mathbb{Z}_{2}\;ou\;x\in \mathbb{Z}_{2}^{\ast },%
			\end{cases}%
		\end{equation}%
		where $x_{t}$ is the coefficient of $2^{t}$ in the $2$-adic expansion of $x$.
	\end{enumerate}
\end{proposition}

\begin{proposition} \cite{ko} \label{pr3}
	\textbf{(Gauss-Legendre multiplication)} 
	\noindent \newline
	For an integer $N$ such as $N\geq 2$ and $(N,p)=1$, let 
	\begin{equation*}
		g_{N}(x)=\Gamma _{q}(x)\prod_{i=0}^{N-1}\Gamma _{q}\left( \frac{i}{N}\right)
		\prod_{i=0}^{N-1}\Gamma _{q}\left( \frac{x+i}{N}\right) ^{-1}.
	\end{equation*}%
	Then $g_{N}(x)=N^{R_{t}(x)-1}N^{(q-1)xR_{t}^{^{\prime }}(x)}$, with $%
	R_{t}^{^{\prime }}(x)=\frac{(x-R_{t}(x))}{q}$.
\end{proposition}

\begin{remark}
	\cite{ko} Just like for the function $\Gamma _{p}$, we have 
	\begin{equation*}
		\prod_{j=1}^{N-1}\Gamma _{q}\left( \frac{j}{N}\right) =%
		\begin{cases}
			\pm 1 & if \quad N odd, \\ 
			\pm 1, or \pm i & if \quad N even.%
		\end{cases}%
	\end{equation*}%
	where $i\in \mathbb{Z}_{p}$, satisfied $i^{2}=-1$.
\end{remark}

\begin{proposition}
	\cite{ko} \label{pr4} Let $a$ and $b$ be integers satisfying $0<a<b<q-1$. 
	For $r\geq 0$,  we define 
	\begin{equation*}
		n_{r}=b\frac{q^{r}-1}{q-1}\quad and \quad m_{r}=a\frac{q^{r}-1}{q-1}.
	\end{equation*}%
	Then for $r>0$ 
	\begin{align}
		\dfrac{\displaystyle\binom{n_{r}}{m_{r}}}{\displaystyle\binom{n_{r-1}}{%
				m_{r-1}}}& =(-1)^{t}(-p)^{v_{p}\left( \binom{b}{a}\right) }\frac{\Gamma
			_{q}(1+n_{r})}{\Gamma _{q}(1+m_{r})\Gamma _{q}(1+n_{r}-m_{r})} \\
		& =(-p)^{v_{p}\left( \binom{b}{a}\right) }\frac{\Gamma _{q}(-m_{r})\Gamma
			_{q}(m_{r}-n_{r})}{\Gamma _{q}(-n_{r})}.
	\end{align}
\end{proposition}
%


In the following, we give our results concerning the combinatorics properties of the generalized $p$-adic gamma function and the coefficients of its Mahler expansion:

\begin{proposition}
	\label{pr4.2.2} The function $\Gamma _{q}$ satisfies the following properties
	
	\begin{enumerate}
		\item For $n$ a positive integer, we have 
		\begin{equation}
			\Gamma _{q}(n+1)=\frac{(-1)^{A_{n}}p^{-B_{n}}n!}{\left[ \frac{n}{q}\right]
				!p^{t\left[ \frac{n}{q}\right] }}.\label{3.30}
		\end{equation}
		with $A_{n}$ and $B_{n}$ as in the proposition \ref{pr2}
		
		\item For $n=mq+\lambda $, such that $\lambda \in \{0,1,\cdots ,q-1\}$ and $%
		m\in \mathbb{N}$, we have 
		\begin{equation}
			\Gamma _{q}(mq+\lambda +1)=\frac{(-1)^{A_{m,\lambda }} p^{-B_{m,\lambda }} (mq+\lambda )!}{%
				m!p^{tm}}, \label{3.31}
		\end{equation}%
		with  
		\begin{equation}
			A_{mq+\lambda }=t+\lambda +mp^{t}+v_{p}(\lambda !)\;  and \;B_{mq+\lambda
			}=m\left( \frac{p^{t}-1}{p-1}\right) -tm-t\left[ \frac{\lambda }{p^{t}}%
			\right] +v_{p}(\lambda !).
		\end{equation}
		
		\item For $n=q^{s}-1$ and $s\in \mathbb{N}$ we have 
		\begin{equation}
			\Gamma _{q}(q^{s})=\frac{(-1)^{A_{q^{s}-1}}p^{-B_{q^{s}-1}}(q^{s}-1)!}{%
				(q^{s-1})!p^{tq^{s-1}}},
		\end{equation}%
		such that 
		\begin{equation*}
			A_{q^{s}-1}=t-1+\sum_{i=0}^{t-1}p^{ts-i}\; and \;B_{q^{s}-1}=%
			\sum_{i=1}^{t}p^{ts-i}-tp^{t(s-1)}.
		\end{equation*}
		
		\item For $n=q^{s}$ and $s\in \mathbb{N}$ we have 
		\begin{equation*}
			\Gamma _{q}(q^{s}+1)=\frac{(-1)^{A_{q^{s}}}p^{-B_{q^{s}}}(q^{s})!}{%
				(q^{s-1})!p^{tq^{s-1}}},
		\end{equation*}%
		such that 
		\begin{equation*}
			A_{q^{s}}-1=A_{q^{s}-1}\;  and \;B_{q^{s}}=B_{q^{s}-1}.
		\end{equation*}
	\end{enumerate}
\end{proposition}

\begin{proof}
	\noindent
	\begin{enumerate}
		\item We know that 
		\begin{equation*}
			n!=\prod_{j=1}^{n}j=\prod_{\substack{ j=1  \\ q\nmid j}}^{n}j\prod 
			_{\substack{ j=1  \\ q\mid j}}^{n}j = \frac{n}{q} \prod_{\substack{ j=1  \\ q\mid j}}^{n}j,
		\end{equation*}%
		and 
		\begin{equation*}
			card\{1\leq j\leq n\;/\;q\mid j\}=\left[ \frac{n}{q}\right] =\left[ \frac{n}{%
				p^{t}}\right] .
		\end{equation*}%
		Thus 
		\begin{equation*}
			\prod_{\substack{ j=1  \\ q\mid j}}^{n}j=\prod_{j=1}^{\left[ \frac{n}{p^{t}}%
				\right] }(p^{t}j)=\left[ \frac{n}{p^{t}}\right] !p^{t\left[ \frac{n}{p^{t}}%
				\right] },
		\end{equation*}%
		therefore 
		\begin{equation*}
			\Gamma _{q}(n+1)=\frac{(-1)^{A_{n}}p^{-B_{n}}n!}{\left[ \frac{n}{q}\right]
				!p^{t\left[ \frac{n}{q}\right] }}.
		\end{equation*}
		
		\item For $n=mq+\lambda $ such that $\lambda \in \{0,1,\cdots ,q-1\}$,
		we have $$\left[ \frac{n}{q}\right] =\left[ \frac{mq+\lambda }{q}\right] =m+%
		\left[ \frac{\lambda }{q}\right] =m$$
		Then 
		\begin{equation*}
			\Gamma _{q}(mq+\lambda +1)=\frac{(-1)^{A_{mq+\lambda }}p^{-B_{mq+\lambda
					}}(mq+\lambda )!}{m!p^{tm}},
			\end{equation*}%
			with 
			\begin{align*}
				A_{mq+\lambda }& =t+\sum_{i=0}^{t-1}\left[ \frac{mq+\lambda }{p^{i}}\right] 
				\\
				& =t+\lambda +mp^{t}+m\sum_{i=1}^{t-1}\left( \frac{1}{p}\right)
				^{i}+\sum_{i=1}^{t-1}\left[ \frac{\lambda }{p^{i}}\right]  \\
				& =t+\lambda +mp^{t}+m\left( \frac{p^{t}-1}{p-1}\right) +\sum_{i=1}^{t-1}%
				\left[ \frac{\lambda }{p^{i}}\right] 
			\end{align*}%
			On other hand we have $v_{p}(\lambda !)=\sum\limits_{i=1}^{t-1}\left[ \frac{\lambda }{%
				p^{i}}\right]$, so 
			\begin{align*}
				(-1)^{A_{mq+\lambda }}& = (-1)^{t+\lambda +mp^{t}+v_{p}(\lambda !)},
			\end{align*}%
			therefore 
			\begin{equation*}
				A'_{m,\lambda }=t+\lambda +mp^{t}+v_{p}(\lambda !),
			\end{equation*}
			
			For $B'_{m,\lambda }$ , we have 
			\begin{align*}
				B_{mq+\lambda }& =\sum_{i=1}^{t}\left[ \frac{mq+\lambda }{p^{i}}\right] -t%
				\left[ \frac{mq+\lambda }{p^{t}}\right]  \\
				& =m\sum_{i=1}^{t-1}\left( \frac{1}{p}\right) ^{i}-tm+\sum_{i=1}^{t}\left[ 
				\frac{\lambda }{p^{i}}\right] -t\left[ \frac{\lambda }{p^{t}}\right]  \\
				& =m\left( \frac{p^{t}-1}{p-1}\right) -tm+\sum_{i=1}^{t}\left[ \frac{\lambda 
				}{p^{i}}\right] -t\left[ \frac{\lambda }{p^{t}}\right],
			\end{align*}%
			then 
			\begin{equation*}
				B'_{m,\lambda }=m\left( \frac{p^{t}-1}{p-1}\right) -tm-t\left[ \frac{\lambda 
				}{p^{t}}\right] +v_{p}(\lambda !).
			\end{equation*}
			
			\item For $n=q^{s}-1$ we have: $\left[ \frac{n}{q}\right] =\left[ \frac{%
				q^{s}-1}{q}\right] =q^{s-1}$, then 
			\begin{equation*}
				\Gamma _{q}(q^{s})=\frac{(-1)^{A_{q^{s}-1}}p^{-B_{q^{s}-1}}(q^{s}-1)!}{%
					(q^{s-1})!p^{tq^{s-1}}},
			\end{equation*}%
			such that 
			\begin{align*}
				A_{q^{s}-1} & = t+\sum_{i=0}^{t-1}\left[ \frac{p^{ts}-1}{p^{i}}\right]  =t-1+\sum_{i=0}^{t-1}p^{ts-i} =t-1+p^{t(s-1)+1}\left( \frac{p^{t}-1}{p-1}\right) ,
			\end{align*}%
			and%
			\begin{align*}
				B_{q^{s}-1} & =\sum_{i=1}^{t}\left[ \frac{q^{s}-1}{p^{i}}\right] -t\left[ \frac{%
					q^{s}-1}{p^{t}}\right] =\sum_{i=1}^{t}p^{ts-i}-tp^{t(s-1)} =p^{t(s-1)}\left( \frac{p^{t}-1}{p-1}\right) -tp^{t(s-1)}.
			\end{align*}
			
			\item For $n=q^{s}$, we have $\left[ \frac{n}{q}\right] =\left[ \frac{q^{s}%
			}{q}\right] =q^{s-1}$, then 
			\begin{equation*}
				\Gamma _{q}(q^{s}+1)=\frac{(-1)^{A_{q^{s}}}p^{-B_{q^{s}}}(q^{s})!}{%
					(q^{s-1})!p^{tq^{s-1}}},
			\end{equation*}%
			such that 
			\begin{align*}
				A_{q^{s}} & =t+\sum_{i=0}^{t-1}\left[ \frac{p^{ts}}{p^{i}}\right]  =t+\sum_{i=0}^{t-1}p^{ts-i} =A_{q^{s-1}}+1,
			\end{align*}%
			and 
			\begin{align*}
				B_{q^{s}} & =\sum_{i=1}^{t}\left[ \frac{p^{ts}}{p^{i}}\right] -t%
				\left[ \frac{p^{ts}}{p^{t}}\right] =\sum_{i=1}^{t}p^{ts-i}-tp^{t(s-1)} =B_{q^{s-1}}.
			\end{align*}
		\end{enumerate}
	
	\end{proof}

	\begin{proposition}
		\label{pr1} \noindent
		
		\begin{enumerate}
			\item We have $\Gamma _{q}(0)=1$ and $|\Gamma _{q}(x)|_{p}=1$ ,  $\forall x\in \mathbb{Z}_{p}$.			
			\item For $p \geq 3$, let  $x,y\in \mathbb{Z}_{p}$, we have \\
			If $|x-y|_{p}=1$, then $|\Gamma _{q}(x)-\Gamma _{q}(y)|_{p}\leq |x-y|_{p}.$ \\
			If $|x-y|_{p}=\frac{1}{p^{s}}$ with $s\geq t$, then  $|\Gamma _{q}(x)-\Gamma _{q}(y)|_{p} \geq |x-y|_{p}$.
		\end{enumerate}
	\end{proposition}
	
	\begin{proof}
		\noindent
		
		\begin{enumerate}
			\item By the definition of the function $\Gamma _{q}$, we have $\Gamma _{q}(0)=\prod_{\ell =0}^{t-1}\Gamma _{p}(h_{\ell }(-1)+1)$,
			
			On the other hand
			\begin{align*}
				-1 & = \sum_{k=0}^{l-1}(p-1)p^{k}+p^{\ell }\left( \sum_{i=0}^{+\infty
				}(p-1)p^{i}\right) ,
			\end{align*}%
			therefore $h_{\ell }(-1)= -1$, so $\Gamma _{q}(0)=\prod_{\ell=0}^{t-1}\Gamma _{p}(0)=1$.
			
			For the second, we again apply the definition of $\Gamma _{q}$ and the fact that $\Gamma _{p}(x)|_{p}=1$, we get 
			\begin{equation*}
				|\Gamma _{q}(x+1)|_{p}=\prod_{\ell =0}^{t-1}|\Gamma _{p}(h_{\ell}(x)+1)|_{p}=1.
			\end{equation*}
			
			\item First, we'll prove the property for positive integers, and then we'll deduce it for p-adic integers by passing to the limit. Let's assume $m,n\in \mathbb{N}$ 
			
			\begin{enumerate}
				\item[\textbf{-}] If $|m-n|_{p}=1$ then 
				\begin{equation*}
					|\Gamma _{q}(m)-\Gamma _{q}(n)|\leq \max (|\Gamma _{q}(m)|_{p},|\Gamma
					_{q}(n)|_{p})\leq 1=|m-n|_{p}
				\end{equation*}
				
				\item[\textbf{-}] If $|m-n|_{p}=\frac{1}{p^{s}}<1$ with $s\in \mathbb{N}$, we have $%
				\exists \mu \in \mathbb{N}$ such that $m=n+\mu p^{s}$ and $(\mu ,p)=1$, therefore 
				\begin{equation*}
					\Gamma _{q}(m)=\Gamma _{q}(n+\mu p^{s})=(-1)^{A}p^{-B}\prod_{\substack{ j=1
							\\ q\nmid j}}^{n+\mu p^{s}-1}j=(-1)^{A}p^{-B}\prod_{\substack{ j=1 \\ q\nmid
							j}}^{n-1}j\prod_{\substack{ j=n \\ q\nmid j}}^{n+\mu p^{s}-1}j 
				\end{equation*}%
				For $s\geq t$ we have $A=A_{1}+A_{2}$ and $B=B_{1}+B_{2}$, with 
				\begin{equation*}
					A_{1}=t+\sum_{i=0}^{t-1}\left[ \frac{n}{p^{i}}\right] , \quad  
					A_{2}=\sum_{i=0}^{t-1}{\mu p^{s-i}},
				\end{equation*}%
				and 
				\begin{equation*}
					B_{1}=\sum_{i=1}^{t}\left[ \frac{n}{p^{i}}\right] -t\left[ \frac{n}{p^{t}}%
					\right] , \quad  B_{2}=\sum_{i=1}^{t}{\mu p^{s-i}}-t{\mu p^{s-t}}
				\end{equation*}%
				Therefore 
				\begin{align*}
					\Gamma _{q}(m)& =\left( (-1)^{A_{1}}p^{-B_{1}}\prod_{\substack{ j=1 \\ %
							q\nmid j}}^{n-1}j\right) \;\left( (-1)^{A_{2}}p^{-B_{2}}\prod_{\substack{ j=n
							\\ q\nmid j}}^{n+\mu p^{s}-1}j\right) \\ 
                            & =\Gamma _{q}(n)\left( (-1)^{A_{2}}p^{-B_{2}}\prod_{\substack{ j=n \\ %
							q\nmid j}}^{n+\mu p^{s}-1}j\right)
				\end{align*}%
				hence 
				\begin{align*}
					|\Gamma _{q}(m)-\Gamma _{q}(n)|_{p}& =|\Gamma _{q}(n)|_{p}\left\vert \left(
					(-1)^{A_{2}}p^{-B_{2}}\prod_{\substack{ j=n \\ q\nmid j}}^{n+\mu
						p^{s}-1}j\right) -1\right\vert _{p}
				\end{align*}%
				Now we rewrite the following product  	
				\begin{equation*}
					\prod_{\substack{ j=n \\ q\nmid j}}^{n+\mu p^{s}-1}j=\prod_{\substack{ j=n
							\\ q\nmid j}}^{n+p^{s}-1}j\;\prod_{\substack{ j=n+p^{s} \\ q\nmid j}}%
					^{n+2p^{s}-1}j\cdots \prod_{\substack{ j=n+(\mu -1)p^{s} \\ q\nmid j}}%
					^{n+\mu p^{s}-1}j
				\end{equation*}%
				So, by the Theorem \ref{3.3}, we have 
				\begin{equation*}
					\prod_{\substack{ j=n+(k -1)p^{s} \\ q\nmid j}}^{n+k p^{s}-1}j \qquad \forall k \geq 1
				\end{equation*}%
				therefore 
				\begin{equation*}
					\prod_{\substack{ j=n \\ q\nmid j}}^{n+\mu p^{s}-1}j\equiv
					(-1)(-1)(-1)\cdots (\mod p^{s})\equiv (-1)^{\mu }(\mod p^{s})
				\end{equation*}%
				hence 
				\begin{equation*}
					(-1)^{A_{2}}p^{-B_{2}}\prod_{\substack{ j=n \\ q\nmid j}}^{n+\mu
						p^{s}-1}j\equiv (-1)^{A_{2}+\mu}p^{-B_{2}}(\mod p^{s})
				\end{equation*}%
				such that 
				\begin{equation*}
					A_{2}+\mu  =\sum_{i=0}^{t-1}\frac{\mu p^{s}}{p^{i}}+\mu  
					 =\mu p^{s+1}\left( \frac{p^{t}-1}{(p-1)p^{t}}\right) +\mu 
				\end{equation*}%
				\begin{equation*}
					B_{2} =\sum_{i=1}^{t}\frac{\mu p^{s}}{p^{i}}-t\frac{\mu p^{s}}{p^{t}} 
					 =\mu p^{s}\left( \frac{p^{t}-1}{(p-1)p^{t}}-\frac{t}{p^{t}}\right).
				\end{equation*}%
				
				So, we have 
				\begin{equation*}
					(-1)^{A_{2}}p^{-B_{2}}\prod_{\substack{ j=n \\ q\nmid j}}^{n+\mu
						p^{s}-1}j\equiv p^{-B_{2}}(\mod p^{s}),
				\end{equation*}%
				then
				\begin{equation*}
					\left\vert (-1)^{A_{2}}p^{-B_{2}}\prod_{\substack{ j=n \\ q\nmid j}}^{n+\mu
						p^{s}-1}j-1\right\vert _{p}=p^{B_{2}} \geq \frac{1}{p^{s}}=|m-n|.
				\end{equation*}
			\end{enumerate}
			
			it means that the property is true for the elements of $\mathbb{Z}_{p}$ by taking the
			limit.
			
		\end{enumerate}
		
		This completes the proof.
	\end{proof}

	\begin{corollary}
		\label{corl4.2.5} For $p\neq 2$, we have 
		\begin{equation*}
			\Gamma _{q}\left( \frac{1}{2}\right) ^{2}=(-1)^{t}.
		\end{equation*}
	\end{corollary}
	
	\begin{proof}
		According to the relationship (\ref{3.34}), for $p\neq 2$ and $x=\frac{1}{2}$, we have 
		\begin{equation*}
			\Gamma _{q}\left( \frac{1}{2}\right) \Gamma _{q}\left( 1-\frac{1}{2}\right)
			=\Gamma _{q}\left( \frac{1}{2}\right) ^{2}=(-1)^{t-1+R_{t}\left( \frac{1}{2}%
				\right) },
		\end{equation*}%
		with 
		\begin{equation*}
			R_{t}\left( \frac{1}{2}\right) \equiv \frac{1}{2}(\mod{q}),
		\end{equation*}%
		such that 
		\begin{equation*}
			q+1\equiv 1(\mod{q})\quad \Longrightarrow \quad \frac{q+1}{2}\equiv \frac{1}{%
				2}(\mod{q}),
		\end{equation*}%
		hence 
		\begin{equation*}
			R_{t}\left( \frac{1}{2}\right) =\frac{q+1}{2}
		\end{equation*}%
		Therefore 
		\begin{equation*}
			\Gamma _{q}\left( \frac{1}{2}\right) ^{2}=(-1)^{t-1+\left( \frac{q+1}{2}%
				\right)} =(-1)^{t}.
		\end{equation*}%
	\end{proof}
	
	
	\begin{proposition} \textbf{(The Mahler expansion)}
		\label{4.2.8} \noindent \newline
		Let the Mahler expansion $\Gamma _{q}(x+1)=\sum_{\eta \geq 0}a_{\eta }\binom{x}{k}$. So
		the coefficients $(a_{\eta })_{\eta \in \mathbb{N}}$ satisfy the relation 
		\begin{equation*}
			\sum_{\eta \geq 0}(-1)^{\eta +t}a_{\eta }\frac{x^{\eta }}{\eta !}=\exp
			_{p}\left( \frac{x^{q}}{q}+x\right) \;\delta _{q}
		\end{equation*}%
		with 
		\begin{equation*}
			\delta _{q}=\sum_{\lambda =0}^{p^{t}-1}(-x)^{\lambda
			}(-1)^{t}(-1)^{v_{p}(\lambda !)}\;p^{v_{p}(\lambda !)-t\left[ \frac{\lambda 
			}{p^{t}}\right] }.
	\end{equation*}
\end{proposition}

\begin{proof}
	We say that the coefficients of the Mahler expansion $(a_{\eta })$ satisfy 
	\begin{equation*}
		\sum_{\eta \geq 0}a_{\eta }\frac{x^{\eta }}{\eta !}=e^{-x}\sum_{\eta \geq
			0}\Gamma _{q}(\eta +1)\frac{x^{\eta }}{\eta !},
	\end{equation*}%
	The aim is to find an expression for the generating series of $(\Gamma _{q}(n+1))_{n\in \mathbb{N}}$  
	\begin{equation*}
		f(x)=\sum_{\eta \geq 0}\Gamma _{q}(\eta +1)\frac{x^{\eta }}{\eta !}.
	\end{equation*}%
	
	However, we start with the Euclidean division of $\eta $ by $q=p^{t}$ is given by
	\begin{equation*}
		\eta =mq+\lambda ,\quad sach that :\;\lambda \in \{0,1,\cdots ,q-1\},
	\end{equation*}%
	Therefore, we have the following cases 
	\begin{align*}
		\lambda =0& \Longrightarrow \eta =mq \\
		\lambda =1& \Longrightarrow \eta =mq+1 \\
		\lambda =2& \Longrightarrow \eta =mq+2 \\
		& \vdots  \\
		\lambda =q-1& \Longrightarrow \eta =mq+(q-1).
	\end{align*}%
	Now, we sum term by term according to the previous cases 
	\begin{equation*}
		f(x)=\sum_{m=0}^{\infty }\frac{\Gamma _{q}(mq+1)}{(mq)!}x^{mq}+\sum_{m=0}^{%
			\infty }\frac{\Gamma _{q}(mq+2)}{(mq+1)!}x^{mq+1}+\cdots +\sum_{m=0}^{\infty
		}\frac{\Gamma _{q}(mq+q)}{(mq+(q-1))!}x^{mq+(q-1)},
	\end{equation*}%
	then 
	\begin{equation*}
		f(x)=\sum_{j=0}^{q-1}\sum_{m=0}^{+\infty }\frac{\Gamma _{q}(mq+j+1)}{%
			(mp^{t}+j)!}x^{mq+j}.
	\end{equation*}%
	We use the relation (\ref{3.30}) and the relation (\ref{3.31}) for: $n=mq+\lambda $, such that $%
	\lambda \in \{0,1,\cdots ,q-1\}$ \newline
	we have $\left[ \frac{n}{q}\right] =\left[ \frac{mq+\lambda }{q}\right] =m+%
	\left[ \frac{\lambda }{q}\right] =m$. Then, for $0\leq \lambda \leq q-1$ 
	\begin{align*}
		\Gamma _{q}(mq+\lambda +1)& =\frac{(-1)^{A_{mq+\lambda }}p^{-B_{mq+\lambda
				}}(mq+\lambda )!}{m!p^{tm}} \\
			& =\frac{(-1)^{t+\lambda +mp^{t}+v_{p}(\lambda !)}p^{m\left( \frac{p^{t}-1}{%
						p-1}\right) -tm+v_{p}(\lambda !)-t\left[ \frac{\lambda }{p^{t}}\right]
				}(mq+\lambda )!}{m!p^{tm}}
		\end{align*}%
		we obtain 
		\begin{align*}
			f(x)& =\sum_{\lambda =0}^{q-1}\sum_{m=0}^{+\infty }\frac{\Gamma
				_{q}(mq+\lambda +1)}{(mq+\lambda )!}x^{mq+\lambda } \\
			& =\sum_{\lambda =0}^{p^{t}-1}\sum_{m=0}^{+\infty }\frac{(-1)^{A_{mq+\lambda
					}}\;(-1)^{mp^{t}-mp^{t}}\;(-1)^{\lambda -\lambda
				}\;(-1)^{t-t}\;p^{-B_{mq+\lambda }}\;(mp^{t}+\lambda )!}{m!\;p^{tm}%
				\;(mp^{t}+\lambda )!}x^{mp^{t}+\lambda } \\
			& =\sum_{\lambda =0}^{p^{t}-1}\sum_{m=0}^{+\infty }\left( \left( \frac{%
				(-1)^{p^{t}}x^{p^{t}}}{p^{t}}\right) ^{m}\frac{1}{m!}\right) (-x)^{\lambda
			}(-1)^{t}(-1)^{v_{p}(\lambda !)}\;p^{v_{p}(\lambda !)-t\left[ \frac{\lambda 
			}{p^{t}}\right] } \\
		& =\exp _{p}\left( \frac{(-x)^{p^{t}}}{p^{t}}\right) (-1)^{t}\;\sum_{\lambda
			=0}^{p^{t}-1}(-x)^{\lambda }(-1)^{t}(-1)^{v_{p}(\lambda
			!)}\;p^{v_{p}(\lambda !)-t\left[ \frac{\lambda }{p^{t}}\right] } \\
		& =\exp _{p}\left( \frac{(-x)^{q}}{q}\right) (-1)^{t}\sum_{\lambda
			=0}^{p^{t}-1}(-x)^{\lambda }(-1)^{t}(-1)^{v_{p}(\lambda
			!)}\;p^{v_{p}(\lambda !)-t\left[ \frac{\lambda }{p^{t}}\right] }
	\end{align*}%
	Now, we put
	\begin{equation*}
		\delta _{q}=\sum_{\lambda =0}^{p^{t}-1}(-x)^{\lambda
		}(-1)^{t}(-1)^{v_{p}(\lambda !)}\;p^{v_{p}(\lambda !)-t\left[ \frac{\lambda 
		}{p^{t}}\right] }
\end{equation*}%
hence 
\begin{equation*}
	e^{-x}f(x)=\exp _{p}\left( \frac{(-x)^{q}}{q}-x\right) (-1)^{t}\;\delta _{q},
\end{equation*}%
then 
\begin{equation*}
	\sum_{\eta \geq 0}a_{\eta }\frac{x^{\eta }}{\eta !}=\exp _{p}\left( \frac{%
		(-x)^{q}}{q}-x\right) (-1)^{t}\;\delta _{q},
\end{equation*}%
we replace $(-x)$ by $x$, it follows 
\begin{equation*}
	\sum_{\eta \geq 0}(-1)^{\eta }a_{\eta }\frac{x^{\eta }}{\eta !}=\exp
	_{p}\left( \frac{x^{q}}{q}+x\right) (-1)^{t}\;\delta _{q},
\end{equation*}%
finally 
\begin{equation*}
	\sum_{\eta \geq 0}(-1)^{\eta +t}a_{\eta }\frac{x^{\eta }}{\eta !}=\exp
	_{p}\left( \frac{x^{q}}{q}+x\right) \;\delta _{q}.
\end{equation*}%
This completes the proof.
\end{proof}

\section*{Acknowledgements} 

This paper is supported by the Scientific Research Project N° C00L03UN180120180006, at the University of Jijel.

Rafik Belhadef\\
LMPA, University of Jijel, BP 98. Jijel, Algeria\\
Corresponding author. Email: belhadefrafik@univ-jijel.dz, rbelhadef@gmail.com
Nour Elhouda SAHALI\\
Department of Mathematics, University of Jijel, BP 98. Jijel, Algeria\\
hodasahali2000@gmail.com
\label{lastpage}

\end{document}